\documentclass[12pt,reqno]{amsart}

\usepackage[usenames]{color}
\usepackage{amssymb}
\usepackage{amsmath}
\usepackage{amsthm}
\usepackage{amsfonts}
\usepackage{amscd}
\usepackage{mathrsfs}
\usepackage{hyperref}
\usepackage{enumerate}
\DeclareMathOperator{\lcm}{lcm}

\theoremstyle{plain}
\newtheorem{theorem}{Theorem}

\newtheorem{lemma}[theorem]{Lemma}
\newtheorem{proposition}[theorem]{Proposition}

\theoremstyle{definition}
\newtheorem{definition}[theorem]{Definition}
\newtheorem{example}[theorem]{Example}

\theoremstyle{remark}

\newcommand\fix{{\rm{Fix}}}

\newcommand\divides\mid
\newcommand{\exactlydivides}{\parallel} 
\newcommand\smalldivides{\mathrel{\kern-2pt\kern3.5pt|}}
\newcommand\notdivides{\mathrel{\kern-3pt\not\!\kern4.3pt\bigm|}}
\newcommand\smallnotdivides{\mathrel{\kern-2pt\not\!\kern3.5pt|}}
\renewcommand{\le}{\leqslant}
\renewcommand{\ge}{\geqslant}
\newcommand\ord{{\rm{ord}}}
\newcommand\reg{{\rm{reg}}}

\let\emptyset\varnothing

\begin{document}

\title[Generating all time-changes]{Generating all time-changes\\
preserving dynamical zeta functions}

\author{Sawian Jaidee}
\address{Department of Mathematics, Faculty of Science, Khon Kaen University, Thailand}
\email{jsawia@kku.ac.th}
\thanks{The named first author thanks Khon Kaen University
for funding a research visit
to Durham University in~2023.}

\author{Jakub Byszewski}
\address{Institute of Mathematics,
Jagiellonian University, Krakow, Poland}
\email{jakub.byszewski@gmail.com}

\author{Thomas Ward}
\address{Department of Mathematical Sciences,
Durham University, Durham, England}
\email{tbward@gmail.com}

\thanks{To the memory of Patrick Moss, who passed
away on the 17th of March 2024}


\date{\today}

\begin{abstract}
We construct a set of topological
generators for the monoid of time-changes
preserving the space of dynamical zeta
functions.
\end{abstract}

\maketitle

\section{Introduction}

A `dynamical system' for our purposes is simply
a map~$T\colon X\to X$ on a set with the property
that
\[
\fix_n(T)=\vert\{x\in X\mid T^nx=x\}\vert<\infty
\]
for all~$n\ge1$, where~$T^n$ denotes the~$n$th
iterate of~$T$.
Associated to any such system is the formal
`dynamical zeta function',
\[
\zeta_T(z)=\exp\sum_{n\ge1}\tfrac{\fix_n(T)}{n}z^n,
\]
and we denote by~$\mathcal{Z}$ the set of all dynamical
zeta functions.
In many settings the variable~$z$ is taken to
lie in~$\mathbb{C}$
and (for example) Tauberian arguments relate analytic properties
of~$\zeta_T$ to orbit-growth properties of~$T$;
for our purposes~$\zeta_T$ is best thought of as
a formal power series in~$z$.
This set enjoys some algebraic properties:
\begin{itemize}
\item $\mathcal{Z}$ is closed under multiplication of functions
in the usual sense
(corresponding to
disjoint unions of the underlying dynamical systems); cryptically
we may write~$\zeta_T\zeta_S=\zeta_{T\sqcup S}$;
\item $\mathcal{Z}$ is closed
under a form of Hadamard multiplication (corresponding to
Cartesian products of the underlying systems), which we
may write as
\[
(\zeta_T*_{\rm{H}}\zeta_S)(z)
=
\exp\sum_{n\ge1}\tfrac{\fix_n(T)\fix_n(S)}{n}z^n
=
\zeta_{T\times S}(z);
\]
\item $\mathcal{Z}$ is not closed under addition or composition
(see Example~\ref{examplenotclosed}).
\end{itemize}
The set~$\mathcal{Z}$ has special properties
because integer sequences of the form~$(\fix_n(T))_{n\ge1}$
(the `realizable' sequences)
are far from arbitrary:
\begin{itemize}
\item $\fix_n(T)\in\mathbb{N}_0$ for all~$n\ge1$;
\item the sequence~$(\fix_n(T))_{n\ge1}$ obeys
an infinite family of
congruence conditions (for example,
the difference~$\fix_2(T)-\fix_1(T)$
is the count of points on the closed orbits of length~$2$
for~$T$
and hence is even);
\item if~$\zeta_T$ is a rational function
then the poles and zeros can only occur
at reciprocals of algebraic integers
(see~\cite[Prop.~1]{MR0271401}).
\end{itemize}

Here we describe a topological
generating set for the monoid of
`symmetries' or `time changes'
of~$\mathcal{Z}$ in the following sense.
The second and third authors together with Patrick
Moss introduced the following definition in~\cite{MR4002553}:
Given a zeta function~$\zeta_T$ and a map~$h\colon\mathbb{N}\to\mathbb{N}$
we define the time-change of~$\zeta_T$ by~$h$
to be the formal function
\[
(h*\zeta_T)(z)
=
\exp\sum_{n\ge1}\tfrac{\fix_{h(n)}{(T)}}{n}z^n.
\]
That is, if iteration of~$T$ is thought of
as the passage of time, then~$h$ changes time for~$T$,
replacing~$T^n$ with~$T^{h(n)}$.
The function~$h*\zeta_T$ may or may not be a
dynamical zeta function, and we define
\[
\mathscr{P}(X,T)=\{h\colon\mathbb{N}\to\mathbb{N}\mid
h*\zeta_{T}\in\mathcal{Z}\}
\]
to be the allowed time-changes for
the dynamical zeta function~$\zeta_T$
(or, equivalently, for the map~$T$).
If~$h\in\mathscr{P}(X,T)$ then there is, by definition,
another dynamical system~$T^{(h)}\colon X^{(h)}\to X^{(h)}$
with
\[
\fix_n(T^{(h)})=\vert\{x\in X^{(h)}\mid
T^{h(n)}x=x\}\vert
\]
for all~$n\ge1$.

Some maps~$h\colon\mathbb{N}\to\mathbb{N}$ may have
this preserving property for all dynamical systems, giving
rise to a collection of maps
\[
\mathscr{P}=\bigcap_{\{(X,T)\}}\mathscr{P}(X,T)
\]
where the intersection is taken over all dynamical systems.
What can be said immediately is that~$\mathscr{P}$ is not
empty, because it contains the map~$h(n)=n$ for all~$n\ge1$.
It is also clear that it is closed under composition, and
so defines a monoid with identity.
The two main results in~\cite{MR4002553} are
as follows:
\begin{itemize}
\item $\mathscr{P}$ is large in that it is uncountable, and
\item $\mathscr{P}$ is small in that the only polynomials it contains are monomials~$n\mapsto an^b$ for~$a,b\in\mathbb{N}_0$.
\end{itemize}

Notice that the presence of monomials means that~$\mathscr{P}$
is not commutative. The maps~$n\mapsto2n^2$
and~$n\mapsto3n$ are in~$\mathscr{P}$ and
do not commute with each other under
composition, for example.

Here we wish to describe a family of
topological generators
of~$\mathscr{P}$ in the following sense.
There is a natural topology on~$\mathscr{P}$
and we will exhibit families
of maps~$\{g_{p,t}\mid p\in\mathbb{P},t\ge0\}$
and~$\{h_{p,t}\mid p\in\mathbb{P},t\ge0\}$
that together generate a dense sub-monoid of~$\mathscr{P}$,
where~$\mathbb{P}$ is the set of rational primes.

\begin{example}\label{examplenotclosed}
The full shift~$T$ on~$2$ symbols
has~$\zeta_T(z)=\frac{1}{1-2z}$.
The Taylor
expansion of~$\log\left(\frac{2}{1-2z}\right)$
has some non-integer coefficient, showing
that~$\mathcal{Z}$ is not closed under addition.
Similarly,~$\log\left(\frac{1}{1-2(\frac{1}{1-2z})}\right)$
has some non-integer coefficient, showing
that~$\mathcal{Z}$ is (unsurprisingly) not closed under composition.
\end{example}

\section{Notation}

The arguments below are elementary
but intricate and potentially
confusing. We will write composition of maps
as
\[
fg(n)=(f\circ g)(n)=f(g(n)),
\]
and carry the formal product notation
as far as writing~$\prod_{i\in I}f_i$
for a composition of functions indexed
by an ordered finite set~$I$.
We write~$\mathbb{P}$ for the
set of prime numbers
written
in their natural order as~$\mathbb{P}=\{p_1,p_2,p_3,\dots\}$.
For a prime~$p\in\mathbb{P}$
write~$p^k\exactlydivides n$ for `exactly divides', meaning
that~$p^k\divides n$ but~$p^{k+1}\notdivides n$.
We write~$\ord_p(n)$ for the maximal power of~$p$
dividing~$n$, so by definition~$p^{\ord_p(n)}\exactlydivides n$,
and we call~$p^{\ord_p(n)}$ the~$p$-part of~$n$.

\section{The topology on~$\mathscr{P}$}

A map~$h\in\mathbb{N}^{\mathbb{N}}$ lies in~$\mathscr{P}$
if and only if it preserves two properties on sequences of non-negative integers~$(a_n)_{n\ge1}$:
If~$\sum_{d\divides n}\mu(n/d)a_d$
is non-negative and divisible by~$n$ for all~$n\ge1$
then~$\sum_{d\divides n}\mu(n/d)a_{h(d)}$
is non-negative (property~(S))
and divisible by~$n$ (property~(D)) for all~$n\ge1$.
We refer to the recent survey~\cite{MR4332826}
for more on this and the
context in which the two
properties---the `Dold congruence'~(D)
and the `sign condition'~(S)---sit.
The space~$\mathbb{N}^{\mathbb{N}}$ of
functions~$\mathbb{N}\to\mathbb{N}$ inherits
a natural topology of pointwise convergence,
meaning that~$h_1$ and~$h_2$ are close in~$\mathscr{P}$
if~$\min\{k\ge1\mid h_1(k)\neq h_2(k)\}$ is large,
and the Dold and sign conditions are clearly
closed, so~$\mathscr{P}$ is a closed sub-monoid
of~$\mathbb{N}^{\mathbb{N}}$.

In order to illustrate the finite
combinatorial nature of the ideas we will use
and the meaning of the definition of~$\mathscr{P}$,
we also give a direct proof that~$\mathbb{N}^{\mathbb{N}}
\setminus\mathscr{P}$
is open.
If~$f\colon\mathbb{N}\to\mathbb{N}$ is
not in~$\mathscr{P}$, this means there is
a dynamical system~$(X,T)$ that
witnesses this fact via
the statement~$f\notin\mathscr{P}(X,T)$.
That is, there is some dynamical system~$(X,T)$
for which~$f*\zeta_{T}$ is not a
dynamical zeta function.
That means~(D) and~(S) hold for
the sequence~$(\fix_T(n))_{n\ge1}$
but there is some~$k\in\mathbb{N}$ such
that
\begin{align*}
\sum_{d\smalldivides k}
\fix_{T}(f(d))\mu(k/d)&\not\equiv0\pmod{k}
\intertext{or}
\sum_{d\smalldivides k}
\fix_{T}(f(d))\mu(k/d)&<0.
\end{align*}
If~$h\in\mathbb{N}^{\mathbb{N}}$ is sufficiently
close to~$f$, then~$h(n)=f(n)$ for all~$n\le k$,
and so the statement above means
that~$(X,T)$ is also a witness to the
statement~$h\in\mathbb{N}^{\mathbb{N}}
\setminus\mathscr{P}$, showing that the
complement of~$\mathscr{P}$ is open.

\section{The maps $g_{p,t}$ and $h_{p,t}$}

In his thesis~\cite{pm} Patrick Moss constructed several
elements of~$\mathscr{P}$, including:
\begin{itemize}
\item the map~$g_p$ defined by~$g_p(n)=\begin{cases}pn
&\mbox{if }p\divides n\\
n&\mbox{if not}\end{cases}$
for any~$p\in\mathbb{P}$; and
\item the map~$n\mapsto n^2$.
\end{itemize}
A natural question is to ask if the
maps~$\{g_p\mid p\in\mathbb{P}\}$
together generate a dense subset of~$\mathscr{P}$,
but this is clearly false. For example,
for any constant~$k\in\mathbb{N}$ the
function~$n\mapsto k$ lies in~$\mathscr{P}$,
and such a map cannot be approximated by
compositions of maps of the form~$g_p$.
The same holds for another almost trivial
observation, namely that the map~$n\mapsto2n$
lies in~$\mathscr{P}$
(this holds simply because~$T^{2n}=(T^2)^n$ for all~$n\ge1$).
Indeed it turns out that much finer control
of how a function changes
the primes dividing~$n$ is essential
to approximate elements of~$\mathscr{P}$.

\begin{definition}\label{definitiongptandhptmaps}
For~$p\in\mathbb{P}$ and~$t\ge0$ define maps~$g_{p,t}$
and~$h_{p,t}$ by
\begin{align*}
g_{p,t}(n)&=\begin{cases}
pn&\mbox{if }\ord_p(n)=t;\\
n&\mbox{if }\ord_p(n)\neq t,
\end{cases}
\intertext{and}
h_{p,t}(n)&=\begin{cases}
n&\mbox{if }\ord_p(n)<t;\\
np^{t-\ord_p(n)}&\mbox{if }\ord_p(n)\ge t.
\end{cases}
\end{align*}
\end{definition}

Notice that
\begin{align}\label{equationAddedExplanation}
\ord_p(h_{p,t}(n))
&=
\min\{t,\ord_p(n)\}
\intertext{and}
\ord_p(g_{p,t}(n))
&\neq
t\label{equationAddedExplanation2}
\end{align}
for any~$n\ge1$.

\begin{lemma}\label{lemmahptandgptinP}
For any~$p\in\mathbb{P}$
and~$t\ge0$ we have~$g_{p,t}\in\mathscr{P}$
and~$h_{p,t}\in\mathscr{P}$.
\end{lemma}

In order to prove this we first assemble
two elementary reductions in the criteria
for membership in~$\mathscr{P}$.

\begin{lemma}\label{lemmasingleorbit}
A map~$f\colon\mathbb{N}\to\mathbb{N}$ lies
in~$\mathscr{P}$ if and only
if~$f*\zeta_T\in\mathcal{Z}$ for
any dynamical system~$T\colon X\to X$
comprising a single orbit.
\end{lemma}

\begin{proof}
A single orbit means that~$(\fix_n(T))=(\reg_k(n))$
for some~$k\ge1$, where
\[
\reg_k(n)=\begin{cases}k&\mbox{if }k\divides n;\\
0&\mbox{if not}\end{cases}
\]
for all~$n\ge1$.
If~$f\in\mathscr{P}$ then~$(\reg_k(f(n))$ is
realizable.
For the reverse direction notice that
properties~(D) and~(S) are additive,
and the periodic orbits of any dynamical
system simply comprise a disjoint union
of single orbits.
\end{proof}

Lemma~\ref{lemmasingleorbit} gives a convenient
way to test candidate functions for membership in~$\mathscr{P}$.

\begin{example}
The function~$h\colon n\to n^n$
time-changes~$\reg_8$ to the sequence
beginning~$(0,0,0,8,0,8,\dots)$ since~$8\divides 4^4$
and~$8\divides 6^6$ and~$8\notdivides k^k$
for~$k=1,2,3,5$. If there were a dynamical
system~$(X,T)$ with this sequence of periodic
point counts, it would have~$2$ orbits of
length~$4$ and no shorter orbits, so~$6\notdivides\fix_T(6)$
is impossible.
\end{example}

\begin{lemma}\label{Pdefinedusingpreimages}
A map~$f\colon\mathbb{N}\to\mathbb{N}$ lies
in~$\mathscr{P}$ if and only
if for any~$k\in\mathbb{N}$
the pre-image
\[
f^{-1}(k\mathbb{N})=\{n\in\mathbb{N}\mid f(n)\in k\mathbb{N}\}
=
\begin{cases}\emptyset&\mbox{or}\\
d_f(k)\mathbb{N}
\end{cases}
\]
for some~$d_f(k)\divides k$.
\end{lemma}

\begin{proof}
By Lemma~\ref{lemmasingleorbit} it is enough to
check the property implies
that the sequence~$(\reg_k(f(n)))$ is realizable
for any~$k\ge1$.
This sequence takes values in~$\{0,k\}$
and so could be identically zero.
The first non-zero term if one exists
must be~$k$, and must correspond to~$\frac{k}{d(k)}$
closed orbit of length~$d(k)$ for some~$d(k)\divides k$.
That is,
\[
\reg_k(f(n))
=
\frac{k}{d(k)}\reg_{d(k)}(n)
\]
for all~$n\ge1$.
It follows that~$f^{-1}(k\mathbb{N})=\emptyset$
or~$f^{-1}(k\mathbb{N})=d(k)\mathbb{N}$ for some~$d(k)\divides
k$ as required.
\end{proof}

\begin{proof}[Proof of Lemma~\ref{lemmahptandgptinP}.]
Let~$k>0$ and~$t\ge1$ be integers.
We claim that
\begin{align*}
h_{p,t}^{-1}(k \mathbb{N})&=\begin{cases}
 \emptyset&\mbox{if }\ord_p(k)\ge t+1;\\
k\mathbb{N}&\mbox{if }\ord_p(k)\le t.
\end{cases}
\end{align*}
Let~$\ord_p(k)\ge t+1$ and suppose that there exists~$n\ge1$
such that
\[
h(n)=kn_1
\]
for some~$n_1\ge1$.
If~$\ord_p(n)<t$
then
\[
\ord_p(n)=\ord_p(k)+\ord_p(n_1)\ge t+1,
\]
a contradiction. If~$\ord_p(n)\ge t$ then,
since~$\ord_p(h(n))=t$
by~\eqref{equationAddedExplanation},
we have
\[
t=\ord_p(kn_1)\ge t+1,
\]
a contradiction. Thus~$h_{p,t}^{-1}(k \mathbb{N})=\emptyset$ if~$\ord_p(k)\ge t+1$.

Now assume that~$\ord_p(k)\le t$
and~$n\ge1$ has~$h_{p,t}(n)=kn_2$ for some~$n_2\ge1$.
If~$\ord_p(n)< t$,
then~$n=h(n)=kn_2 \in k\mathbb{N}$.
If~$\ord_p(n)\ge t$, then
\[
h_{p,t}(n)=np^{t-\ord_p(n)}=
kn_2
\]
so~$n=kn_2p^{\ord_p(n)-t} \in k\mathbb{N}.$
Hence~$h_{p,t}^{-1}(k \mathbb{N})\subseteq  k\mathbb{N}.$

If~$h_{p,t}^{-1}(k\mathbb{N})\neq k\mathbb{N}$
then there is some~$m=kn\in k\mathbb{N}$
 with~$h_{p,t}(m)\notin k\mathbb{N}$.
If~$\ord_{p}(kn)<t$ then~$h_{p,t}(kn)=kn\in k\mathbb{N}$,
so we must have~$\ord_p(kn)\ge t$.
Then
\[
\ord_p(h_{p,t}(kn))=\min\{t, \ord_p(kn)\}=t<\ord_p(k)
\]
since~$h_{p,t}(kn)\notin k\mathbb{N}$ means
that the~$p$-part of~$h_{p,t}(kn)$ is a strict
divisor of the~$p$-part of~$k$.
This contradicts~$\ord_p(k)\le t$,
so~$h_{p,t}^{-1}(k \mathbb{N})=k\mathbb{N}.$

We claim that
\begin{align*}
g_{p,t}^{-1}(k \mathbb{N})&=\begin{cases}
k\mathbb{N}&\mbox{if }\ord_p(k)=t
\mbox{ or }\ord_p(k)<t\mbox{ or }\ord_p(k)>t+1;\\
\frac{k}{p}\mathbb{N}&\mbox{if }\ord_p(k)=t+1.
\end{cases}
\end{align*}

Assume first that~$\ord_p(k)=t$
and~$g_{p,t}(n)=kn_1$.
If~$\ord_p(n_1)=1$ then~$n=\frac{kn_1}{p}$
or~$n=kn_1$ so~$n\in k\mathbb{N}$.
By~\eqref{equationAddedExplanation2}
we have~$\ord_p(n_1)\neq 0$,
and if~$\ord_p(n_1)>1$ then~$n=kn_1\in k\mathbb{N}$.
Hence~$g_{p,t}^{-1}(k\mathbb{N})\subset k\mathbb{N}$.
Clearly~$g_{p,t}(k\mathbb{N})\subset k\mathbb{N}$,
so we have~$g_{p,t}^{-2}(k\mathbb{N})=k\mathbb{N}$
as claimed.

If~$\ord_p(k)=t+1$ and~$g_{p,t}(n)\in k\mathbb{N}$
then~$pn\in k\mathbb{N}$ or~$n\in k\mathbb{N}$,
so~$g_{p,t}^{-1}(k\mathbb{N})\subset\frac{k}{p}\mathbb{N}$.
On the other hand~$g_{p,t}(k\mathbb{N})=k\mathbb{N}\subset
\frac{k}{p}\mathbb{N}$, so
\[
g_{p,t}^{-1}(k\mathbb{N})=\frac{k}{p}\mathbb{N}.
\]

If~$\ord_p(k)<t$ or~$\ord_p(k)>t+1$
then
and~$g_{p,t}(n)=kn_1$ implies
that either
\[
\ord_p(n)+1=\ord_p(k)+\ord_p(n_1)
\]
or
\[
\ord_p(n)=\ord_p(k)+\ord_p(n_1).
\]
In either case we deduce that~$\ord_p(n)\neq t$,
so~$g_{p,t}(n)=n$ and
therefore~$g_{p,t}(k\mathbb{N})=k\mathbb{N}$.
\end{proof}


In order to see how these maps together generate
other elements in~$\mathscr{P}$ we go through
the four examples we have already seen, namely~$g_p$,
the squaring map, a constant map, and the
doubling map.
Notice that each of these maps can be
written uniquely using prime factorizations in the form
\[
\prod_{p\in\mathbb{P}}p^{v_p}\longmapsto
\prod_{p\in\mathbb{P}}p^{d_p(v_p)}
\]
for certain maps~$d_p\colon\mathbb{N}_0\to\mathbb{N}_0$.

\begin{example}\label{examplegp}
The map~$g_p$ has~$d_p(v_p)=v_p+1$ and~$d_q(v_q)=v_q$ for all
primes~$q\neq p$.
The composition
\[
g_{p,1}g_{p,2}\cdots g_{p,k}(n)
\]
multiplies~$n$ by~$p$ if~$\ord_p(n)\le k$,
and so it converges to~$g_p$ as~$k\to\infty$.
\end{example}

\begin{example}\label{examplesquaring}
The map~$n\mapsto n^2$ lies in~$\mathscr{P}$,
and has~$d_p(v_p)=2v_p$ for all~$p\in\mathbb{P}$.
To see how this works, notice first that~$n\mapsto g_{p,1}(n)$
squares the~$p$-part of~$n$ if~$\ord_p(n)\le1$,
and~$n\mapsto g_{p,1}g_{p,3}g_{p,2}(n)$ squares
the~$p$-part of~$n$ if~$\ord_p(n)\le2$, and so on.
Thus
\[
f_p(n)
=
\underbrace{g_{p,1}}_{\times p\text{ if }\ord_p=1;}
\underbrace{g_{p,3}g_{p,2}}_{\times p^2\text{ if }\ord_p=2;}
\underbrace{g_{p,5}g_{p,4}g_{p,3}}_{\times p^3\text{ if }\ord_p=3;}
\cdots
\underbrace{g_{p,2k-1}\cdots g_{p,k+1}g_{p,k}}_{\times p^k\text{ if }\ord_p=k}(n)
\]
squares the~$p$-part of~$n$ if~$\ord_p(n)\le k$.
So the composition
\[
f_{p_k}f_{p_{k-1}}\cdots f_{2}(n)
\]
squares the~$p$-part of~$n$ if~$\ord_p(n)\le k$
for all primes~$p\le p_k$.
That is, this composition agrees with~$n\mapsto n^2$
for
\[
n\in\{p_1^{j_1}\cdots p_k^{j_k}\mid
j_i\le k\mbox{ for all }1\le i\le k\}.
\]
Thus~$n\mapsto n^2$ lies in the closure of
the sub-monoid generated by the maps~$g_{p,t}$
and~$h_{p,t}$.
\end{example}

\begin{example}\label{exampleconstant2}
The map~$n\mapsto2$ has~$d_2(v_2)=1$ and~$d_p(v_p)=0$ for all
primes~$p\neq2$.
We claim that for~$k\ge3$ we have
\[
h_{p_k,0}h_{p_{k-1},0}\cdots h_{7,0}h_{5,0}h_{3,0}h_{2,1}g_{2,0}
(n)=2
\]
for all~$n<p_k$.
To see the claim,
notice that the first map~$g_{2,0}$ multiplies by~$2$
if~$n$ is odd, resulting in an even number.
Then~$h_{2,1}$ applied
to any even number~$n$ gives~$n2^{1-\ord_p(n)}$,
a number divisible by~$2$ exactly once.
For any prime~$p<p_k$ the map~$h_{p,0}$
removes all powers of~$p$ that can be removed.
It follows that~$n\mapsto 2$ lies in the closure
of the sub-monoid generated by the maps~$g_{p,t}$
and~$h_{p,t}$.
\end{example}

\begin{example}\label{exampledoubling}
The map~$n\mapsto2n$ lies in~$\mathscr{P}$
(it corresponds to replacing~$T\colon X\to X$
with~$T^2\colon X\to X$).
This map has
\[
d_p(v_p)
=
\begin{cases}
v_2+1&\mbox{if }p=2;\\
v_p&\mbox{if }p\neq2.
\end{cases}
\]
To approximate this map, we need to
multiply by~$2$ exactly once, and then
make no further changes.
Since~$n\mapsto g_{2,t}(n)$ multiplies
by~$2$ only when~$\ord_2(n)=t$, we have
\[
g_{2,0}g_{2,1}g_{2,2}\cdots g_{2,k}(n)=2n
\]
for all~$n<2^{k+1}$. It follows that~$n\mapsto2n$
lies in the closure of the sub-monoid generated
by the maps~$g_{p,t}$.
\end{example}

The maps~$g_{p,t}$ and~$h_{p,t}$ have
somewhat messy orbits on a given dynamical
zeta function, but explicit
calculations are possible for the simpler
map~$g_{p}$ for~$p\in\mathbb{P}$.
Starting with a rational zeta
function this produces algebraic dynamical
zeta functions.

\begin{example}
The function~$\zeta_{T}(z)=\frac{1}{1-2z}$ is
the dynamical zeta function of the full shift
on~$2$ symbols.
For a fixed~$p\in\mathbb{P}$ we have
\begin{align*}
\zeta_{T^{(g_p)}}(z)
&=
\exp\left(
\sum_{k=1}^{\infty}\frac{z^{pk}}{pk}2^k
+
\sum_{n=1}^{\infty}\frac{z^n}{n}2^n
-
\sum_{k=1}^{\infty}\frac{z^{pk}}{pk}2^{pk}
\right)\\
&=
\left(\frac{1}{1-2z^p}\right)^{1/p}\left(\frac{1}{1-2z}\right)
\left(1-2^pz^p\right).
\end{align*}
\end{example}

\section{Commutation relations}

As mentioned earlier, the monoid~$\mathcal{P}$ is
certainly not commutative, and this is
reflected in properties of the maps~$g_{p,t}$ and~$h_{p,t}$.

\begin{lemma}\label{lemmacommutationrelations}
For~$p,p_1,p_2\in\mathbb{P}$ and~$t,t_1,t_2\ge0$
the following relations hold.
\begin{enumerate}
\item[{\rm{(a)}}] $g_{p_1,t_1}$ and~$h_{p_2,t_2}$
commute
if~$p_1\neq p_2$.
\item[{\rm{(b)}}]
$g_{p_1,t_1}$ and~$h_{p_2,t_2}$
commute
if~$p_1=p_2$ and~$t_1\neq t_2$.
\item[{\rm(c)}] $g_{p_1,t_1}$ and~$g_{p_2,t_2}$ commute
if~$p_1\neq p_2$.
\item[{\rm(d)}] $h_{p_1,t_1}$ and~$h_{p_2,t_2}$
commute.
\item[{\rm{(e)}}] $g_{p,t}$ and~$h_{p,t}$ do not commute
but satisfy the relation
\[
h_{p,t+1}(g_{p,t}(n))=g_{p,t}(h_{p,t}(n)).
\]
\end{enumerate}
\end{lemma}

\begin{proof}
Referring to Definition~\ref{definitiongptandhptmaps},
we may consider four cases:
\begin{align}
\ord_{p_1}(n)&=t_1,\label{case1}\\
\ord_{p_1}(n)&\neq t_1,\label{case2}\\
\ord_{p_2}(n)&<t_2,\label{case3}\\
\ord_{p_2}(n)&\ge t_2\label{case4}.
\end{align}

For~(a) we have~$p_1\neq p_2$
and compute the composition in each of the
four possible cases as follows:
\begin{center}
\begin{tabular}{ |c|c|c| }
\hline
 case & $g_{p_1,t_1}(h_{p_2,t_2}(n))$ & $h_{p_2,t_2}(g_{p_1,t_1}(n))$ \\
 \hline
 \eqref{case1}\&\eqref{case3} & $p_1n$ & $p_1n$ \\
 \hline
 \eqref{case1}\&\eqref{case4} & $np_1p_2^{t_2-\ord_{p_2}(n)}$ & $np_1p_2^{t_2-\ord_{p_2}(n)}$ \\
 \hline
 \eqref{case2}\&\eqref{case3} & $n$ & $n$\\
 \hline
 \eqref{case2}\&\eqref{case4} & $np_2^{t_2-\ord_{p_2}(n)}$ & $np_2^{t_2-\ord_{p_2}(n)}$\\
 \hline
\end{tabular}
\end{center}

For~(b)
some additional steps are required. For purposes
of brevity we suppress the single prime~$p$ from the notation,
writing~$g_{t_1}$ and~$h_{t_2}$.
The formula for~$h_{t_2}(g_{t_1}(n))$ falls into
four cases as follows.
\begin{enumerate}
\item[$(\alpha)$] If~$\ord_p(n)=t_1$ and~$t_1+1<t_2$
then~$h_{t_2}(g_{t_1}(n))=pn$.
\item[$(\beta)$] If~$\ord_p(n)=t_1$ and~$t_1+1\ge t_2$
then
\[
h_{t_2}(g_{t_1}(n))=pnp^{t_2-\ord_p(pn)}=
np^{t_2-\ord_p(n)}.
\]
\item[$(\gamma)$] If~$\ord_{p}(n)\neq t_1$ and~$\ord_p(n)<t_2$ then~$h_{t_2}(g_{t_1}(n))=n$.
\item[$(\delta)$] If~$\ord_{p}(n)\neq t_1$ and~$\ord_p(n)\ge t_2$ then~$h_{t_2}(g_{t_1}(n))=np^{t_2-\ord_p(n)}$.
\end{enumerate}
The formula for~$g_{t_1}(h_{t_2}(n))$ falls into
three cases as follows.
\begin{enumerate}
\item[$(\epsilon)$] If~$\ord_p(n)<t_2$ and~$\ord_p(n)=t_1$
then~$g_{t_1}(h_{t_2}(n))=pn$.
\item[$(\zeta)$] If~$\ord_p(n)<t_2$ and~$\ord_p(n)\neq t_1$
then~$g_{t_1}(h_{t_2}(n))=n$.
\item[$(\eta)$] If~$\ord_p(n)\ge t_2$ and~$t_2\neq t_1$
then~$g_{t_1}(h_{t_2}(n))=np^{t_2-\ord_p(n)}$.
\end{enumerate}
Notice that the fourth possibility
($\ord_p(n)<t_1$ and~$t_2=t_1$) is excluded
by the hypothesis that~$t_1\neq t_2$.
To prove~(b) we have to verify that all the
cases agree.
\begin{itemize}
\item In~$(\epsilon)$ we have~$t_1<t_2$. If~$t_1+1<t_2$
then we match~$(\alpha)$. If~$t_1+1\ge t_2$
then~$t_1+1=t_2$ since~$t_1<t_2$,
so by~$(\beta)$ we have
\[
h_{t_2}(g_{t_1}(n))
=
np^{t_2-\ord_p(n)}
=
np^{t_1+1-\ord_p(n)}=np
\]
which matches up.
\item The cases~$(\zeta)$ and~$(\gamma)$ coincide.
\item Assume we are in case~$(\eta)$.
If~$\ord_p(n)=t_1$ then~$t_1\ge t_2$ so we have~$t_1+1\ge t_2$
and we match~$(\beta)$.
\item If~$\ord_p(n)\neq t_1$ then we match~$(\delta)$.
\end{itemize}

The statement in~(c) is clear.

For~(d), the statement is clear if~$p_1\neq p_2$
since~$h_{p,t}(n)$ only involves the~$p$-part
of~$n$. If~$p_1=p_2=p$ then~$h_{p,t}=h_t$ only
influences the~$p$-part, so it is enough to
calculate~$\ord_p$ as follows.
If
\[
\ord_p(n)\ge t_2\ge t_1
\]
then~$\ord_p(h_{t_1}(h_{t_2}(n)))=t_1$.
If~$\ord_p(n)\ge t_2$ and~$t_2<t_1$,
then
\[
\ord_p(h_{t_1}(h_{t_2}(n)))=t_2.
\]
If~$\ord_p(n)<t_2$ and~$\ord_p(n)\ge t_1$ then
\[
\ord_p(h_{t_1}(h_{t_2}(n)))=t_1.
\]
If~$\ord_p(n)<t_2$ and~$\ord_p(n)< t_1$ then
\[
\ord_p(h_{t_1}(h_{t_2}(n)))=\ord_p(n).
\]
The same four cases for~$h_{t_2}(h_{t_1}(n))$
match up.

For~(e) we compute directly that
\begin{align*}
h_{t+1}(g_{t}(n))
&=
\begin{cases}
h_{t+1}(pn)&\mbox{if }\ord_p(n)=t;\\
h_{t+1}(n)&\mbox{if }\ord_p(n)\neq t
\end{cases}\\
&=
\begin{cases}
pnp^{t+1-\ord_p(n)}=pn&\mbox{if }\ord_p(n)=t;\\
np^{t+1-\ord_p(n)}&\mbox{if }\ord_p(n)>t;\\
n&\mbox{if }\ord_p(n)<t
\end{cases}
\end{align*}
and
\begin{align*}
g_{t}(h_{t}(n))
&=
\begin{cases}
g_{t}(n)&\mbox{if }\ord_p(n)<t;\\
g_{t}(np^{t-\ord_p(n)})&\mbox{if }\ord_p(n)\ge t
\end{cases}\\
&=
\begin{cases}
n&\mbox{if }\ord_p(n)<t;\\
pnp^{t-\ord_p(n)}=
np^{t+1-\ord_p(n)}&\mbox{if }\ord_p(n)> t;\\
pn&\mbox{if }\ord_p(n)=t
\end{cases}
\end{align*}
as required.
\end{proof}

\section{Characterizing $\mathscr{P}$}

In addition to Lemma~\ref{Pdefinedusingpreimages},
it will be useful to have another characterization
of membership in~$\mathscr{P}$.

\begin{proposition}\label{Pdefinedusingfactorization}
A map~$f\in\mathbb{N}^{\mathbb{N}}$ lies in~$\mathscr{P}$
if and only if it can be written
in the form
\begin{equation}\label{JakubClaim}
f\colon n=\prod_{p\in\mathbb{P}}p^{v_p(n)}
\longmapsto
\prod_{p\in\mathbb{P}}p^{d_p(v_p(n))}
\end{equation}
for functions~$d_p\colon\mathbb{N}_0\to\mathbb{N}_0$
satisfying the following properties
\begin{enumerate}
\item[\rm(i)] $d_p(0)=0$ for all but finitely many~$p\in\mathbb{P}$.
\item[\rm(ii)] $d_p$ is non-decreasing for every~$p\in\mathbb{P}$.
\item[\rm(iii)] If~$d_p$ is unbounded then~$d_p(i)\ge i$ for all~$i\ge0$.
\item[\rm(iv)] If~$d_p$ is bounded, then~$d_p(i)\ge i$ for all~$i\le\max\{d_p(i)\mid i\in\mathbb{N}_0\}$.
\end{enumerate}
\end{proposition}

\begin{proof}
Suppose first that~$f\in\mathscr{P}$.
We claim that
\begin{equation}\label{JakubClaimI}
m\divides n \Longrightarrow f(m)\divides f(n)
\end{equation}
for all~$m,n\in\mathbb{N}$.
Since~$f(m)\in f(m)\mathbb{N}$ the set~$f^{-1}(f(m)\mathbb{N})$ contains~$m$,
so by Lemma~\ref{Pdefinedusingpreimages}
we have
\[
f^{-1}(f(m)\mathbb{N})=d(f(m))\mathbb{N}
\]
for some~$d(f(m))\divides f(m)$. If~$m\divides n$
it follows that~$n\in d(f(m))\mathbb{N}$ since~$m\in d(f(m))\mathbb{N}$, so~$f(n)\in f(m)\mathbb{N}$ as claimed in~\eqref{JakubClaimI}.

We next claim that
\begin{equation}\label{JakubClaimII}
\gcd(m,n)=1
\Longrightarrow
f(mn)=\lcm(f(m),f(n))
\end{equation}
for all~$m,n\in\mathbb{N}$.
Let~$k=\lcm(f(m),f(n))$,
so~$f(m)\divides k$ and~$f(n)\divides k$.
Since~$m\divides mn$ and~$n\divides mn$
we have~$f(m)\divides f(mn)$ and~$f(n)\divides f(mn)$,
so~$k\divides f(mn)$.
Suppose that~$k\neq f(mn)$,
and pick~$p\in\mathbb{P}$ and~$t\ge1$
so that~$p^t\divides f(mn)$ and~$p^t\notdivides k$.
Since~$p^t\divides f(mn)$
we have
\[
f^{-1}(p^t\mathbb{N})=\{\ell\ge1\mid f(\ell)\in p^t\mathbb{N}\}
\ni mn,
\]
so~$f^{-1}(p^t\mathbb{N})=d(p^t)\mathbb{N}$
for some~$d(p^t)\divides p^t$.
Also~$f(mn)\in p^t\mathbb{N}$ so
\[
mn\in f^{-1}(p^t\mathbb{N})
=d(p^t)\mathbb{N},
\]
so~$d(p^t)\divides mn$.
Since~$\gcd(m,n)=1$ it follows that~$d(p^t)\divides m$
or~$d(p^t)\divides n$.
We have~$f(d(p^t))\in p^t\mathbb{N}$,
so if~$d(p^t)\divides m$ then~$p^t\divides f(d(p^t))\divides \divides f(m)$
by~\eqref{JakubClaimI}.
The same argument applies to~$n$,
so~$p^t\divides f(m)$ or~$p^t\divides f(n)$,
which contradicts the assumption that~$p^t\notdivides k=\lcm(f(m),f(n))$ and proves~\eqref{JakubClaimII}.

Next we claim that
\begin{equation}\label{JakubClaimIII}
p\divides\tfrac{f(n)}{f(1)}
\Longrightarrow
p\divides n
\end{equation}
for any~$p\in\mathbb{P}$ and~$n\ge1$.
If~$n=1$ the statement is trivial.
If~\eqref{JakubClaimIII} is known
for all prime powers, then given~$n\in\mathbb{N}$
we may factorize as~$n=\prod_{p\in\mathbb{P}}p^{v_p}$.
By~\eqref{JakubClaimII} it follows that
\[
\tfrac{f(n)}{f(1)}
=
\tfrac{\lcm(f(2^{v_2}),f(3^{v_3}),\dots)}{f(1)},
\]
so if a prime~$q$ divides~$\tfrac{f(n)}{f(1)}$
then~$q\divides\tfrac{f(p^{v_p})}{f(1)}$
for some~$p$ by~\eqref{JakubClaimIII} for prime
powers, so~$q\divides n$.
Assume therefore that~$n=p^t$
for some~$p\in\mathbb{P}$ and~$t\ge1$
and there is a prime~$q\neq p$
such that~$q^s\exactlydivides f(1)$
and~$q^{s+1}\divides f(p^t)$,
so
\[
q\divides\tfrac{f(p^t)}{f(1)}.
\]
Then~$f^{-1}(q^{s+1}\mathbb{N})\ni p^t$
so~$f^{-1}(q^{s+1}\mathbb{N})\neq\emptyset$.
Hence
\[
f^{-1}(q^{s+1}\mathbb{N})=
d(q^{s+1})\mathbb{N}
\]
for some~$d(q^{s+1})\divides q^{s+1}$.
Since~$d(q^{s+1})\divides p^t$ we
must have~$d(q^{s+1})=1$, so~$f(1)\in q^{s+1}\mathbb{N}$,
contradicting the assumption that~$q^s\exactlydivides f(1)$.
This proves~\eqref{JakubClaimIII} for prime
powers and hence for all~$n\ge1$.

The claims~\eqref{JakubClaimI}
and~\eqref{JakubClaimIII} imply that for each~$p\in\mathbb{P}$
there is a non-decreasing function~$d_p\colon\mathbb{N}_0\to\mathbb{N}_0$
with~$f(1)=\prod_{p\in\mathbb{P}}p^{d_p(0)}$
and
\[
f(p^t)=f(1)p^{d_p(t)-d_p(0)},
\]
so in particular
with~$d_p(0)=0$ for all but finitely many~$p$.
By~\eqref{JakubClaimII}, it follows that~$f$
has the form in~\eqref{JakubClaim}.

If for some~$p\in\mathbb{P}$
there are~$i,j\in\mathbb{N}$ with~$d_p(j)\ge i$, then
the set~$f^{-1}(p^i\mathbb{N})$
contains~$p^j$, so~$f^{-1}(p^i\mathbb{N})
=d(p^i)\mathbb{N}$
for some~$d(p^i)$ which must
divide both~$p^i$ and~$p^j$, so~$d(p^i)=p^{\ell}$
for some~$\ell$ with
\[
0\le\ell\le\min\{i,j\}.
\]
It follows that
\[
f(d(p^i))
=
f(p^{\ell})
=
f(1)p^{d_p(\ell)-d_p(0)}\in p^i\mathbb{N},
\]
so~$d_p(j)\ge d_p(\ell)\ge i$.
This shows that~$f$ satisfies~(iii) and~(iv).
\end{proof}

\section{Generating $\mathscr{P}$}

Our main result is in effect a generalization of
Examples~\ref{examplegp}--\ref{exampledoubling}
to find topological generators of~$\mathscr{P}$.

\begin{theorem}
The monoid~$\mathscr{P}'$ generated by
the maps~$g_{p,t}$ and~$h_{p,t}$
for~$p\in\mathbb{P}$
and~$t\ge0$ is dense
in~$\mathscr{P}$.
\end{theorem}

\begin{proof}
Using the commutation relations in
Lemma~\ref{lemmacommutationrelations}
any term~$h_{*}$ can be moved
left in
any word in~$\mathscr{P}'$, and then
the terms can be grouped according to
the prime involved. That is, any element
of~$\mathscr{P}'$ has the form
\begin{equation}\label{equationlongcomposition}
\prod_{p\in\mathbb{S}_1\subset\mathbb{P}}\prod_{*}h_{p,*}
\prod_{p\in\mathbb{S}_2\subset\mathbb{P}}\prod_{*}g_{p,*}
\end{equation}
where~$\prod_{*}$ denotes a composition of finitely many
choices for the parameter~$t$ and~$\mathbb{S}_1,
\mathbb{S}_2$ are finite sets of primes.
Given~$f\in\mathscr{P}$, we may write it in
the form~\eqref{JakubClaim}
from Lemma~\ref{Pdefinedusingfactorization}.

For any~$p\in\mathbb{P}$ with unbounded~$d_p$
and~$t\ge0$ with~$d_p(t)\ge1$
we may use a composition
\begin{equation}\label{equationg-gadget}
\prod_{*}g_{p,*}
=
g_{p,t}g_{p,t+1}\cdots g_{p,d_p(t)-2}g_{p,d_p(t)-1}
\end{equation}
to multiply~$n$ by~$p^{d_p(t)-t}$ if~$\ord_p(n)=t$
as required.
We compose expressions of this
form for~$t=0,1,\dots,t^*$
and for all primes~$p<t^*$ say.

If~$d_p$ is bounded with~$d_p(v_p)=M_p$
for all~$v_p>t'$ then we compose
terms~\eqref{equationg-gadget}
for~$t=0,1,\dots,t'$ and for each
such~$p$.
This is composed with
the map~$h_{p,t'}$
to reduce any higher power of~$p$
to~$p^{t'}$
and with~$h_{q,0}$
for
\[
q\in\{p\in\mathbb{P}\mid
q<p_{t^*},
d_p(t)=0\mbox{ for all }t\ge0\}.
\]

Composing all these maps as
in~\eqref{equationlongcomposition}
produces a map whose
value at~$n$ agrees with~$f$
as long as~$\ord_p(n)\le t^*$ for all~$p<t^*$
and~$n<t^*$.
As~$t^*$ was arbitrary we
can approximate~$f\in\mathscr{P}$ by
elements of~$\mathscr{P}'$.
\end{proof}

A natural question is to ask if this
is precise in the following sense:
Does~$\mathscr{P}'$
have a presentation with generators
from Definition~\ref{definitiongptandhptmaps}
and relations from Lemma~\ref{lemmacommutationrelations}?

\bibliographystyle{plain}
\bibliography{refs}

\begin{thebibliography}{1}

\bibitem{MR0271401}
R.~Bowen and O.~E. Lanford, III.
\newblock Zeta functions of restrictions of the shift transformation.
\newblock In {\em Global {A}nalysis ({P}roc. {S}ympos. {P}ure {M}ath., {V}ols.
  {XIV}, {XV}, {XVI}, {B}erkeley, {C}alif., 1968)}, pages 43--49. Amer. Math.
  Soc., Providence, R.I., 1970.

\bibitem{MR4332826}
J.~Byszewski, G.~Graff, and T.~Ward.
\newblock Dold sequences, periodic points, and dynamics.
\newblock {\em Bull. Lond. Math. Soc.}, 53(5):1263--1298, 2021.

\bibitem{MR4002553}
S.~Jaidee, P.~Moss, and T.~Ward.
\newblock Time-changes preserving zeta functions.
\newblock {\em Proc. Amer. Math. Soc.}, 147(10):4425--4438, 2019.

\bibitem{pm}
P.~Moss.
\newblock {\em The arithmetic of realizable sequences}.
\newblock PhD thesis, University of East Anglia, 2003.

\end{thebibliography}
\end{document}